\documentclass[a4paper,12pt]{article}
\usepackage[top=2.5cm,bottom=2.5cm,left=2.5cm,right=2.5cm]{geometry}
\usepackage{amsfonts}
\usepackage{latexsym,amsmath,amssymb,amsfonts,epsfig,psfrag,url,graphics,ifpdf,multicol}
\usepackage{amsthm}
\usepackage{ifsym}
\usepackage{color,colordvi,amscd,amsthm} %eepic
\usepackage{tikz}
\usetikzlibrary{decorations.pathreplacing}
\usepackage[numbers,sort&compress]{natbib}
\usepackage{indentfirst,graphics,epsfig}
\usepackage{graphicx}
\usepackage{graphics}
\usepackage{graphicx,psfrag}
\usepackage{graphics}
\usepackage{mathrsfs}
\usepackage{tabularx}
\usepackage{booktabs}
\usepackage{floatrow}
\usepackage{multirow}
\usepackage{graphicx}
\usepackage{caption}
\usepackage[ruled]{algorithm2e}

\newfloatcommand{capbtabbox}{table}[][\FBwidth]

\usepackage{etoolbox}

\newtheorem{thm}{Theorem}

\newtheorem{letterthm}{Theorem}

\newcommand{\thmletter}[1]{\ifcase#1\or A\or B\or C\or D\or E\else ?\fi}

\newtheorem{lem}{Lemma}

\newtheorem{cl}{Claim}
\newtheorem{coro}{Corollary}

 \theoremstyle{definition}
\newtheorem{defn}{Definition}

\newtheorem{rem}{Remark}
\newcommand{\setmytheorem}[1]{%
	\ifcase#1\or A\or B\or C\or D\or E\or F\or G\or H\or I\or J\else
	The number is too large for lettered theorems!\fi
	\setcounter{mytheorem}{#1}% ?????
}

\def\pf{\noindent{\bf Proof.\ }}
\def\qed{{\hfill\rule{4pt}{7pt}}}

\allowdisplaybreaks

\usepackage{indentfirst,subfig}

\begin{document}
 %\thispagestyle{empty}
    %\rule{0cm}{0.5mm}
    \captionsetup[figure]{labelfont={bf},name={Fig.},labelsep=period}

\begin{center} {\large Maximal independent sets in graphs with a given matching number}
\end{center}
\pagestyle{empty}

	\begin{center}
	{
		{\small Yongtang Shi$^1$, Jianhua Tu$^{2,}$\footnote{Corresponding author.\\\indent \ \  E-mail: tujh81@163.com (J. Tu)}, Ziyuan Wang$^2$}\\[2mm]
		
		{\small $^1$Center for Combinatorics and LPMC, Nankai University, Tianjin 300071, China\\			
			$^2$School of Mathematics and Statistics, Beijing Technology and Business University, \\
			\hspace*{1pt} Beijing 100048, China}\\[2mm]}
	
\end{center}

\begin{center}
\begin{abstract}
A maximal independent set in a graph $G$ is an independent set that cannot be extended to a larger independent set by adding any vertex from 
$G$. This paper investigates the problem of determining the maximum number of maximal independent sets in terms of the matching number of a graph. We establish the maximum number of maximal independent sets for general graphs, connected graphs, triangle-free graphs, and connected triangle-free graphs with a given matching number, and characterize the extremal graphs achieving these maxima. 

\vskip 3mm


\end{abstract}
%\end{minipage}
\end{center}

%%%%%%%%%%%%%%%%%%%%%%%%%%%%%%%%%%%%%%%%%%%%

\section{Introduction}

An independent set in a graph is a set of vertices that induces a subgraph without any edges. A \textit{maximal independent set} (MIS, for short) is one that cannot be a proper subset of any larger independent set. For a graph $G$, define $MIS(G)$ as the set of all MISs in $G$, and let $mis(G)=|MIS(G)|$.

Let $K_n$, $C_n$, $P_n$, and $K_{1,n-1}$ denote the complete graph, the cycle, the path, and the star on $n$ vertices, respectively. Around 1960, Erd\H{o}s and Moser raised the problem of determining the maximum value of $mis(G)$ in terms of the order of $G$, and the extremal graphs. The well-known result due to Miller and Muller \cite{Miller1960} and Moon and Moser \cite{Moon1965}, which answers this problem, shows that
for any graph of order $n\geq 2$, 
  \[
  mis(G)\leq\begin{cases}
  	3^{n/3},\ &\text{if\ } n\equiv0 \pmod{3};\\
  	4\cdot 3^{(n-4)/3},\ &\text{if\ } n\equiv1 \pmod{3};\\
  	2\cdot 3^{(n-2)/3},\ &\text{if\ } n\equiv2 \pmod{3}.
  \end{cases}
  \]
  
For two graphs $G$ and $H$, we define $G\cup H$ to be their disjoint union and write $kG$ for the disjoint union of $k$ copies of $G$. Then the extremal graphs achieving the maximum value are $\frac{n}{3}K_3$, $\frac{n-4}{3}K_3\cup K_4$ or $\frac{n-4}{3}K_3\cup 2K_2$, and $\frac{n-2}{3}K_3\cup K_2$ in the three respective cases.

Since then, there has been interest in further exploring the maximum number of MISs in specific families of graphs and in identifying the extremal graphs that attain these numbers. Researchers have determined the maximum value of $mis(G)$ for various graph families of a given order, such as trees, forests, connected graphs, bipartite graphs, unicyclic graphs, and graphs with at most $r$ cycles. For these results, we refer to \cite{Furedi1987, Griggs1988, Hujtera1993, Koh2008, Liu1994, Palmer2023, Sagan1988, Sagan2006, Taleskii2022, Wilf1986, Wloch2008, Ying2006}.

A matching of a graph $G$ is a set of edges with no shared endpoints. A maximum matching is one that has the largest possible number of edges among all matchings in the graph. The matching number of $G$, denoted by $\mu(G)$, is the number of edges in a maximum matching in $G$.

 We observe that for many graph families, such as general graphs, connected graphs, and triangle-free graphs, the extremal graphs that attain the maximum number of maximal independent sets among all graphs of order $n$ always have large matching numbers. Motivated by this, we investigate the maximum number of maximal independent sets in terms of the matching number. We determine this number for general graphs, connected graphs, triangle-free graphs, and connected triangle-free graphs. Specifically, we show that  

(1) for a general graph $G$ with matching number $t$, 
\[mis(G)\leq 3^t,\]

(2) for a connected graph $G$ with matching number $t$, 
\[mis(G)\leq \begin{cases}
	3\ &\text{for\ } t=1,\\
	3^{t-1}+2^{t-1}\ &\text{for\ } t\geq 2,
\end{cases}\]

(3) for a triangle-free graph $G$ with matching number $t$, 
\[mis(G)\leq \begin{cases}
	5^{t/2} & \text{for } t \text{ even},\\
	2\cdot 5^{(t-1)/2} & \text{for } t \text{ odd},
\end{cases}\]

(4) for a connected triangle-free graph $G$ with matching number $t$, 
\[mis(G)\leq \begin{cases}
	5 & \text{for } t=2,\\
	2\cdot(5^{(t-2)/2} + 3^{(t-2)/2}) & \text{for } t\geq 4 \text{ even}, \\
	5^{(t-1)/2} + 3^{(t-1)/2} & \text{for } t \text{ odd}.
\end{cases}\]
For the families of graphs under consideration, we also characterize the corresponding extremal graphs that possess the maximum number of MISs.

It should be noted that while Hoang and Trung \cite{Hoang2019} have previously established the maximum number of MISs in general graphs with a given matching number, our approach, distinct from theirs, offers a different perspective.

The structure of the rest of this paper is as follows: In Section \ref{sec2}, we establish the maximum number of MISs in general graphs with a given matching number and identify the extremal graphs that achieve this maximum, as detailed in Theorem \ref{thm1}. Section \ref{sec3} begins with Theorem \ref{thm2}, which outlines the maximum number of MISs in connected graphs with a given matching number and the extremal graphs that reach this number. We then present the necessary preliminaries and proceed to prove Theorem \ref{thm2}. In Sections \ref{sec4} and \ref{sec5}, we tackle the cases of triangle-free graphs and connected triangle-free graphs, determining their maximum number of MISs in these graphs with a given matching number and characterizing the extremal graphs that attain these maxima, as elaborated in Theorems \ref{thm3} and \ref{thm4}, respectively.

%%%%%%%%%%%%%%%%%%%%%%%%%%%%%%%%%%%%%%%%%%%%%%%%%%%%%%%%%

\section{Maximal independent sets in general graphs with a given matching number}\label{sec2}

Let $G$ be a graph. For a vertex $v\in V(G)$, we define the open neighborhood of $v$ as $N_G(v)=\{u\in V(G)| uv\in E(G)\}$ and closed neighborhood of $v$ as $N_G[v]=N_G(v)\cup\{v\}$. The degree of $v$, $d_G(v)$, is the cardinality of $N_G(v)$. If $d_G(v)=1$, $v$ is referred to as a leaf of $G$. A vertex that is not a leaf and is adjacent to a leaf is called a supported vertex. To simplify notation without causing confusion, we may abbreviate $N_G(v)$, $N_G[v]$, and $d_G(v)$ to $N(v)$, $N[v]$, and $d(v)$ respectively. For a subset of vertices $S\subseteq V(G)$, the induced subgraph $G[S]$ consists of the vertices in $S$ and all edges between them in $G$. The graph $G-S$ is the result of removing all vertices in $S$ from $G$. When $S$ consists of a single vertex $v$, we write $G-\{v\}$ as $G-v$ for brevity. 

Let $M$ be a matching in $G$, the vertices incident to an edge of $M$ are said to be {\it saturated} by $M$, while the others are {\it unsaturated}. An $M$-alternating path is a path that alternates between edges that are in $M$ and edges that are not. An $M$-augmenting path is an $M$-alternating path that starts and ends with unsaturated vertices.

The following theorem is a fundamental theorem in matching theory.

\begin{letterthm}\cite{Berge1957}\label{thmA}
	A matching $M$ in a graph $G$ is a maximum matching if and only if $G$ has no $M$-augmenting paths.
\end{letterthm}

For a graph $G$, define $I(G)$ as the set of all {\it independent sets} in $G$, and let $i(G)=|I(G)|$. It is important to note that the empty set is also considered an independent set.

\begin{lem}\label{lem2.1}
Let $G$ be a connected graph with matching number $t$. Then
	\[mis(G)\leq 3^t.\]
Furthermore, when $t\geq 2$, the inequality is strict.
\end{lem}

\pf Let $M$ be a maximum matching in $G$, and let $V(M)$ be the set of vertices of $G$ that are saturated by $M$. The subgraph induced by $V(M)$ in $G$ is simply written by $G_M$. Let $F=V(G)\setminus V(M)$. Since $M$ is a maximum matching in $G$, $F$ is an independent set of $G$.

It is easy to see that 
\[i(G_M)\leq 3^t\]
with equality if and only if no two edges of $M$ are joined by an edge of $G$ (that is, $M$ is an induced matching of $G$).

We show that 
\[mis(G) \leq i(G_M)\] 
by constructing an injective mapping $f: MIS(G)\rightarrow I(G_M)$. Let $I \in MIS(G)$, then we define $f(I) = I \cap V(M)$.
Let $I_1$ and $I_2$ be two MISs in $G$ such that $f(I_1)=f(I_2)$, i.e., $I_1\cap V(M)=I_2\cap V(M)$. For any vertex $v\in F$ and any MIS $I$ in $G$, $v\in I$ if and only if $(I\cap V(M))\cap N(v)=\emptyset$. Thus, $I_1=I_2$, showing that $f$ is injective. Hence, 
\[mis(G)\leq i(G_M)\leq 3^t.\]

Furthermore, suppose that $t\geq 2$. Based on the analysis in the previous paragraph, we can derive an important fact: If $I'$ is an independent set of $G_M$, then there is at most one MIS $I$ in $G$ such that $I\cap V(M)=I'$. 

We distinguish the following two cases.

{\bf Case 1.} $i(G_M)<3^t$ or there is no MIS $I$ in $G$ such that $I\cap V(M)=\emptyset$.

In this case, it is easy to see that \[mis(G)<3^t.\]

{\bf Case 2.} $i(G_M)=3^t$ and there exists a MIS $I$ in $G$ such that $I\cap V(M)=\emptyset$.

In this case, the MIS $I$ is actually the set $F$, and moreover, $M$ is an induced matching of $G$. Let
\[M:=\{u_1v_1,u_2v_2,\cdots,u_tv_t\},\ F:=\{w_1,\cdots, w_k\}.\]

We will prove that for each edge $u_iv_i$ (where $1\leq i\leq t$) of $M$, there exists a vertex in $F$ to which both endpoints $u_i$ and $v_i$ of the edge are uniquely adjacent.

For example, consider the edge $u_1v_1$. Since $G$ is connected, without loss of generality, assume that $u_1$ is adjacent to at least one vertex in $F$, say $w_1$. That is, assume $u_1w_1 \in E(G)$. Firstly, $v_1$ is adjacent to at least one vertex in $F$, otherwise $F\cup\{v_1\}$ would also be an independent set in $G$, contradicting the maximality of $I=F$. Secondly, the vertex adjacent to $v_1$ in $F$ can only be $w_1$, otherwise we can find an $M$-augmenting path in $G$, which is a contradiction. Finally, apart from $w_1$, neither $u_1$ nor $v_1$ have any adjacent vertices in $F$, otherwise we would be able to find an $M$-augmenting path in $G$, leading to a contradiction.

Since $G$ is connected, $F$ contains only one vertex, i.e., $F=\{w_1\}$. Thus, $G$ is a graph on $2t+1$ vertices obtained by connecting the vertex $w_1$ to both endpoints of each edges in $M$. Now, since $t\geq 2$,
\[mis(G)=2^t+1<3^t.\]

The proof of Lemma \ref{lem2.1} is complete. \qed

\begin{thm}\label{thm1}
If $G$ is a graph with matching number $t$, then
\[mis(G)\leq 3^t,\]
with equality if and only if $G\cong tK_3\cup r K_1$ for some natural number $r$.
\end{thm}
\pf Assume that $G$ has a total of $k$ non-trivial components $G_1$, $\cdots$, $G_k$ with $\mu(G_1)=t_1$, $\cdots$, $\mu(G_k)=t_k$, where a non-trivial component is a component with at least two vertices. Clearly, $t_1+t_2+\cdots+t_k=t$. Then, by Lemma \ref{lem2.1},
\[mis(G)=\prod\limits_{i=1}^k mis(G_i)\leq \prod\limits_{i=1}^k 3^{t_i}=3^t.\]
Furthermore, if $mis(G)=3^t$, then $t_1=t_2=\cdots=t_k=1$, which implies that any non-trivial component of $G$ is a star or a $K_3$. For a star $S$, $mis(S)=2$. Thus, $mis(G)=3^t$ if and only if $G\cong tK_3\cup r K_1$ for some natural number $r$.
\qed

\section{Maximal independent sets in connected graphs with a given matching number}\label{sec3}

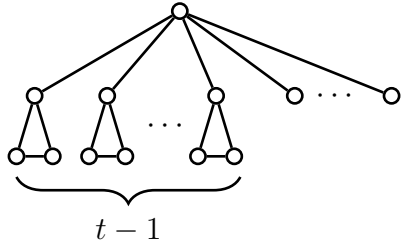
\begin{figure}[H]
	\begin{tikzpicture}
		[line width = 1pt, scale=0.8,
		empty/.style = {circle, draw, fill = white, inner sep=0mm, minimum size=2mm}, full/.style = {circle, draw, fill = black, inner sep=0mm, minimum size=2mm}, full1/.style = {circle, draw, fill = black, inner sep=0mm, minimum size=0.5mm}]
		\node [empty] (a) at (0,2.4) {};
		
		\node [empty] (b) at (-2.4,1) {};
		\node [empty] (c) at (-2.7,0) {};
		\node [empty] (d) at (-2.1,0) {};
		
		\node [empty] (e) at (-1.2,1) {};
		\node [empty] (f) at (-1.5,0) {};
		\node [empty] (g) at (-0.9,0) {};
		
		\node [label=above:$\cdots$] at (-0.2,0) {};
		
		%\node [full1] at (-0.5,1) {};
		%\node [full1] at (-0.2,1) {};
		%\node [full1] at (0.1,1) {};
		\node [empty] (h) at (0.6,1) {};
		\node [empty] (i) at (0.3,0) {};
		\node [empty] (j) at (0.9,0) {};
        \node [empty] (k) at (1.9,1) {};
        %\node [full1] at (2.8,1) {};
		%\node [full1] at (2.4,1) {};
		%\node [full1] at (2.6,1) {};
		
		\node [label=above:$\cdots$] at (2.6,0.5) {};
        \node [empty] (l) at (3.5,1) {};
        
		\draw (a) -- (b);
		\draw (a) -- (e);
		\draw (a) -- (h);
		\draw (b) -- (c);
		\draw (b) -- (d);
		\draw (e) -- (f);
		\draw (e) -- (g);
		\draw (h) -- (i);
		\draw (h) -- (j);
        \draw (c) -- (d);
        \draw (f) -- (g);
        \draw (i) -- (j);
        \draw (a) -- (k);
        \draw (a) -- (l);
        
	\draw[decorate,decoration={brace,mirror,amplitude=10pt},xshift=0pt,yshift=-4pt]
(-2.7,-0.2) -- (1,-0.2) node [black,midway,yshift=-0.7cm] {${t-1}$};

	\end{tikzpicture}
	\caption{An example of a graph in $\mathcal{E}_t$}
	\label{fig1}
\end{figure}

For any integer $t\geq 2$, let $\mathcal{E}_t$ denote the set of connected graphs with matching number $t$, which are constructed by connecting a central vertex to exactly one vertex in each of $t-1$ triangles, and additionally connecting this central vertex to $\ell$ isolated vertices, where $\ell\geq1$. An example of a graph in $\mathcal{E}_t$ is illustrated in Figure \ref{fig1}.

Let $A_5$ be a graph with $5$ vertices, constructed by connecting a central vertex to every vertex in each of $2$ copies of $K_2$, as illustrated in Figure \ref{fig2}.

\begin{figure}[H]
	\begin{tikzpicture}
		[line width = 1pt, scale=0.8,
		empty/.style = {circle, draw, fill = white, inner sep=0mm, minimum size=2mm}, full/.style = {circle, draw, fill = black, inner sep=0mm, minimum size=2mm}, full1/.style = {circle, draw, fill = black, inner sep=0mm, minimum size=0.5mm}]
		\node [empty] (a) at (0,2) {};
		\node [empty] (b) at (2,2) {};
        \node [empty] (c) at (1,1) {};
        \node [empty] (d) at (0,0) {};
        \node [empty] (e) at (2,0) {};
        \draw (a) -- (b);
        \draw (a) -- (c);
        \draw (b) -- (c);
        \draw (c) -- (d);
        \draw (c) -- (e);
        \draw (d) -- (e);
	\end{tikzpicture}
	\caption{$A_5$}
	\label{fig2}
\end{figure}

Let
\begin{align*}
\begin{split}
\mathcal{H}_t:= \left \{
\begin{array}{ll}
    \{K_3\},                  & \text{if}\ t=1, \\
    \{A_5, C_5, K_5\}\cup \mathcal{E}_2,                  & \text{if}\ t=2, \\
    \mathcal{E}_t,    & \text{if}\ t\geq 3.
\end{array}
\right.
\end{split}
\end{align*}

\begin{thm}\label{thm2}
For any connected graph $G$ with matching number $t$, 
\begin{align*}
	\begin{split}
	mis(G)\leq h(t):= \left \{
		\begin{array}{ll}
			3                 & \text{for}\ t=1, \\
			3^{t-1}+2^{t-1}    & \text{for}\ t\geq 2,
		\end{array}
		\right.
	\end{split}
\end{align*}
with equality if and only if $G$ is isomorphic to a graph in $\mathcal{H}_t$.
\end{thm}

\baselineskip=0.24in
%%%%%%%%%%%%%%%%%%%%%%%%%%%%%%%%%%%%%%%%%%%%%%%%%%%

\subsection{Preliminaries}\label{sec3.1}

In the subsection, we provide the preliminaries needed to prove Theorem \ref{thm2}.

\begin{lem}\label{lem3.1}
Let $G$ be a graph. For any induced subgraph $H$ in $G$,
\[mis(H)\leq mis(G).\]
\end{lem}

\pf For any induced subgraph $H$ of $G$, every MIS $I'$ of $H$ can be extended in $G$ to form at least one MIS $I$ such that $I \cap V(H) = I'$. This leads to the inequality 
\[mis(H)\leq mis(G).\]
\qed

\begin{lem}\label{lem3.2}
Let $G$ be a graph with matching number $t$, and let $v$ be a vertex in $G$. Then,

(1) if $v$ is saturated by all maximum matchings of $G$, $\mu(G-v)=t-1$,

 (2) otherwise, $\mu(G-v)=t$.

\end{lem}

\pf Obviously,
\[t-1=\mu(G)-1\leq \mu(G-v)\leq \mu(G)=t.\]

(1) Suppose, to the contrary, that $\mu(G-v)=t$. There exists a maximum matching $M$ in $G-v$ of size $t$. Thus, $M$ is also a maximum matching of $G$ and $v$ is not saturated by $M$, a contradiction. So, $\mu(G-v)=t-1$.

(2) Let $M$ be a maximum matching of $G$ that does not saturate the vertex $v$. Then, $M$ is also a maximum matching of $G-v$. Thus, $\mu(G-v)=t$.

\qed

\begin{lem}\label{lem3.3}
Let $G$ be a graph, and let $v$ be a vertex in $G$ that is saturated by all maximum matchings of $G$. Let $G'$ be the graph obtained from $G$ by adding $k\geq 1$ new vertices and joining these vertices to $v$. Then, $\mu(G')=\mu(G)$. 
\end{lem}
\pf Obviously,
\[\mu(G)\leq \mu(G')\leq \mu(G)+1.\]
Suppose, to the contrary, that $\mu(G')=\mu(G)+1$. By Lemma \ref{lem3.2}, $\mu(G-v)=\mu(G)-1$. Thus, $\mu(G'-v)=\mu(G-v)=\mu(G)-1$. Now, we have $\mu(G')=\mu(G'-v)+2$, a contradiction. So, $\mu(G')=\mu(G)$. 
\qed

Let $H$ be a bipartite graph with bipartition $(X,Y)$. For $S\subseteq X$, define the surplus of $S$ by $\sigma(S)=|N_H(S)|-|S|$. The minimum surplus of non-empty subsets of $X$ is called the surplus of $X$.

Let $G$ be a graph. Define $D(G)$ as the set of vertices that are unsaturated by at least one maximum matching of $G$. Let $A(G)$ be the set of vertices in 
$V(G)$ that are not in $D(G)$ and are adjacent to at least one vertex in $D(G)$. Finally, let $C(G)=V(G)\setminus (A(G)\cup D(G))$. A near-perfect matching in a graph is a matching that saturates all vertices of the graph except for exactly one. If removing any vertex from $G$ results in a graph with a perfect matching, then $G$ is called factor-critical.

\begin{letterthm}\label{thmB}(The Gallai-Edmonds Structure Theorem \cite{Lovasz1986}). For a graph $G$, let $D(G)$, $A(G)$ and $C(G)$ be defined as above. Then:
\begin{itemize}
	\item[(1)] each component of $G[D(G)]$ is factor-critical.
	
	\item[(2)] there is a perfect matching in $G[C(G)]$.
	
	\item[(3)] in the bipartite graph obtained from $G$ by deleting all vertices of $C(G)$, removing all edges in $G[A(G)]$, and contracting each component of $G[D(G)]$ to a single vertex, with bipartition $(X,Y)$ where $X=A(G)$, the surplus of $X$ is positive.

	\item[(4)] every maximum matching in $G$ contains a perfect matching of $G[C(G)]$, a near-perfect matching of each component of $G[D(G)]$ and matches all vertices of $A(G)$ with vertices in distinct components of $G[D(G)]$.	
\end{itemize}
\end{letterthm}

\begin{coro}\label{coro1}
If there is no vertex in a graph $G$ that is saturated by all maximum matchings, then each component of $G$ is factor-critical.
\end{coro}

\pf This corollary can be directly derived from Theorem \ref{thmB}. \qed

\begin{defn}
If $n<6$ let $c(n)=n$; if $n\geq 6$ let 
\begin{align*}
	c(n)=\begin{cases}
		2\cdot 3^{\frac{n-3}{3}}+2^{\frac{n-3}{3}},\ &\text{if}\ n\equiv0\pmod{3},\\
		3^{\frac{n-1}{3}}+2^{\frac{n-4}{3}},\ &\text{if}\ n\equiv1\pmod{3},\\
		4\cdot 3^{\frac{n-5}{3}}+3\cdot 2^{\frac{n-8}{3}},\ &\text{if}\ n\equiv2\pmod{3}.
	\end{cases}
\end{align*}
\end{defn}

\begin{letterthm}\label{thmC}\cite{Griggs1988}
For any connected graph $G$ of order $n$,
\[mis(G)\leq c(n).\]	
\end{letterthm}

\begin{coro}\label{coro2}
For any connected graph $G$ of order $n\geq 11$,
\[mis(G)\leq 3^{\frac{n-1}{3}}+2^{\frac{n-4}{3}}.\]
\end{coro}

\pf This corollary can be directly derived from Theorem \ref{thmC}. \qed

%%%%%%%%%%%%%%%%%%%%%%%%%%%%%%%%%%%%%%%%%%%%%%%

\subsection{Proof of Theorem \ref{thm2}}\label{sec3.2}

\pf Let $g(t)$ be the maximum number of MISs among all connected graphs with matching number $t$.

\begin{rem}
	According to Lemma \ref{lem2.1}, $g(t)\leq 3^t$ for any positive integer $t$, implying that $g(t)$ is a finite integer.
\end{rem}

 For any positive integer $t$ and any graph $G\in \mathcal{H}_t$, by direct calculations, 
\[mis(G)=h(t),\]
which implies that $g(t)\geq h(t)$.

We prove by induction on $t$ that $g(t)\leq h(t)$, which implies that $g(t)=h(t)$.

When $t=1$, because a connected graph with matching number $1$ is either a triangle or a star. A triangle has 3 MISs, while a star has 2 MISs. Therefore, $g(1)=3=h(1)$.

Suppose that $t\geq 2$, and for any connected graph $F$ with matching number $k\leq t-1$, $mis(F)\leq h(k)$. Let $G$ be a connected graph with matching number $t$ and $g(t)$ MISs.

We distinguish the following two cases.

\noindent {\bf Case 1.} There exists a vertex $v$ in $G$ that is saturated by all maximum matchings.

We obtain a graph $G'$ from $G$ by adding a new vertex and joining the new vertex to $v$. By Lemmas \ref{lem3.1} and \ref{lem3.3}, $\mu(G')=\mu(G)=t$ and $mis(G')\geq mis(G)=g(t)$, which implies that $mis(G')=g(t)$.

Let $U$ be the set of all leaves in $G'$ which are adjacent to the vertex $v$. The MISs of $G'$ can be divided into two categories: one includes vertex $v$, and the other includes vertex set $U$. Thus, 
\begin{align*}
g(t)&=mis(G')\\
&=mis(G'-N_{G'}[v])+mis(G'-v-U)\\
&=mis(G'-N_{G'}[v])+mis(G'-v).
\end{align*}

It is easy to see that the vertex $v$ is saturated by all maximum matchings in $G'$. Thus, by Lemma \ref{lem3.2}, $\mu(G'-v)=t-1$. 

When $t=2$, $\mu(G'-v)=1$ and $mis(G'-v)\leq 3$.
Furthermore, if $mis(G'-v)=3$, then $G'-v$ is a disjoint union of a triangle and some isolated vertices, and $mis(G'-N_{G'}[v])\leq 2$ (Since $G'$ is connected), which implies
\[g(2)=mis(G')=mis(G'-N_{G'}[v])+mis(G'-v)\leq 2+3=h(2).\]
Moreover, in Case 1, for $t=2$, $G'$ has $g(2)=h(2)$ MISs if and only if $G'$ is isomorphic to a graph in $\mathcal{E}_2$.

When $t\geq 3$, let $H$ be a component of $G'-v$ that has maximum matching number among all components. Let $\alpha:=mis(H)$ and $\beta:=mis(H-(N_{G'}(v)\cap V(H)))$. Let $G_1:=G'-v-V(H)$ and $G_2:=G'-\left[N_{G'}[v]\cup V(H)\right]$. Clearly, $\mu(H)\geq1$ and $V(G_1)\geq 1$.

We now prove that $\mu(H)=1$. Suppose, to the contrary, that $\mu(H)\geq 2$. Since $\mu(H)\leq 
\mu(G'-v)=t-1$, by inductive hypothesis and Lemma \ref{lem3.1},
\[\beta\leq \alpha\leq 3^{\mu(H)-1}+2^{\mu(H)-1}.\]

We introduce an operation on $G'$ which involves first removing all vertices in $H$ from $G'$, subsequently adding $\mu(H)$ new triangles, and finally connecting the vertex $v$ to exactly one vertex in each of these newly added triangles. The new graph obtained through this process is denoted as $G''$. It is easy to see that $G''$ is a connected graph with matching number $t$. We will prove that $g(t)=mis(G')<mis(G'')$, which contradicts that $g(t)$ is the maximum number of MISs among all connected graphs with matching number $t$. 
Specifically, 
\begin{align*}
	mis(G')=&mis(G'-v)+mis(G'-N_{G'}[v])\\
	=&\alpha\cdot mis(G_1)+\beta\cdot mis(G_2)\\
	\leq &(3^{\mu(H)-1}+2^{\mu(H)-1})\cdot mis(G_1)+(3^{\mu(H)-1}+2^{\mu(H)-1})\cdot mis(G_2)\\
	=& 3^{\mu(H)}\cdot mis(G_1)+2^{\mu(H)}\cdot mis(G_2)\\
	&-\left[(2\cdot 3^{\mu(H)-1}-2^{\mu(H)-1})\cdot mis(G_1)-(3^{\mu(H)-1}-2^{\mu(H)-1})\cdot mis(G_2)\right],
\end{align*}
since $(2\cdot 3^{\mu(H)-1}-2^{\mu(H)-1})>(3^{\mu(H)-1}-2^{\mu(H)-1})\geq 1$, $mis(G_1)\geq 1$ (since $|V(G_1)|\geq1$), and $mis(G_1)\geq mis(G_2)$ (by Lemma \ref{lem3.1}), we have 
\[(2\cdot 3^{\mu(H)-1}-2^{\mu(H)-1})\cdot mis(G_1)-(3^{\mu(H)-1}-2^{\mu(H)-1})\cdot mis(G_2)>0,\]
and
\[mis(G')<3^{\mu(H)}\cdot mis(G_1)+2^{\mu(H)}\cdot mis(G_2)=mis(G'').\]
Hence, we have proved that $\mu(H)=1$.

Since $\mu(H)=1$, each component of $G'-v$ is either a triangle, a star, or an isolated vertex. Thus, 
\begin{align*}
	mis(G')&=mis(G'-v)+mis(G'-N_{G'}[v])\\
	&=\alpha\cdot mis(G_1)+\beta\cdot mis(G_2)\\
	&\leq 3\cdot mis(G_1)+2\cdot mis(G_2).
\end{align*}
Furthermore, the equality in the last inequality holds if and only if $H$ is a triangle and has exactly one vertex adjacent to $v$. Since each component with matching number 1 in $G'-v$ exhibits symmetry, $mis(G')$ achieves its maximum value only when each of these components is a triangle and has exactly one vertex adjacent to $v$. In other words, $mis(G')$ is maximized only when $G'$ is isomorphic to a graph in $\mathcal{E}_t$. Therefore, 
\[g(t)=mis(G')\leq h(t).\]
Moreover, in Case 1, for $t\geq 3$, $G'$ has $g(t)=h(t)$ MISs if and only if $G'$ is isomorphic to a graph in $\mathcal{E}_t.$

%%%%%%%%%%%%%%%%%%%%%%%%%%%%%%%%%%%%
\vskip 3mm

\noindent {\bf Case 2.} There is no vertex in $G$ that is saturated by all maximum matchings.

In this case, by Corollary \ref{coro1}, $G$ is factor-critical. Thus, $|V(G)|=2t+1$.

When $t=2$, by Theorem \ref{thmC}, $g(2)\leq c(5)=5=h(2)$. A few simple calculations show that if a graph $G$ on $5$ vertices is factor-critical and has $5$ MISs, then $G$ is isomorphic to a graph in $\{A_5,K_5,C_5\}.$

When $t=3$, by Theorem \ref{thmC},
\[g(3)\leq c(7)=11<3^2+2^2=h(3).\]  

When $t=4$, by Theorem \ref{thmC},
\[g(4)\leq c(9)=22<3^3+2^3=h(4).\]  

When $t\geq 5$, by Corollary \ref{coro2},
\[g(t)\leq 3^{\frac{2t}{3}}+2^{\frac{2t-3}{3}}<3^{t-1}+2^{t-1}=h(t).\] 

Thus far, we have shown that $g(t)\leq h(t)$, which also demonstrates that $g(t)=h(t)$.

\vskip 3mm

Next, we show that for a connected graph $G$ with matching number $t$, $mis(G)=h(t)$ if and only if $G$ is isomorphic to a graph in $\mathcal{H}_t$.

When $t=1$, it is easy to see that for a connected graph $G$ with matching number $1$, $mis(G)=h(1)=3$ if and only if $G\cong K_3$. 

For any integer $t\geq 2$, let $L_t$ be a connected graph with $3(t-1)+1$ vertices and matching number $t$ constructed as follows: a central vertex is connected to one vertex of each of $t-1$ triangles. One can easily verify that 
\[mis(L)=3^{t-1}<3^{t-1}+2^{t-1}=h(t).\]

When $t=2$, let $G$ be a connected graph with matching number $2$ and $h(2)=5$ MISs. From the details of proving that $g(t)\leq h(t)$, we can observe that $G$ must belong to one of the following cases:

(1) There exists a vertex $v$ in $G$ such that adding a leaf to $v$ results in a graph $G'$ isomorphic to one in $\mathcal{E}_2$.

(2) $G$ is factor-critical. Moreover, $G$ is isomorphic to a graph in $\{A_5,K_5,C_5\}$.

If case (1) occurs, then because $mis(L_2)<h(2)$, we can infer that $G'$, being isomorphic to a graph in $\mathcal{E}_2$, implies that $G$ itself is isomorphic to a graph in $\mathcal{E}_2$.

In conclusion, $G$ is isomorphic to a graph in $\{A_5,K_5,C_5\}\cup\mathcal{E}_2$.

When $t\geq 3$, let $G$ be a connected graph with matching number $t$ and $h(t)$ MISs. From the details of proving that $g(t)\leq h(t)$, we can observe that there exists a vertex $v$ in $G$ such that adding a leaf to $v$ results in a graph $G'$ isomorphic to one in $\mathcal{E}_t$.

Since $mis(L_t)<h(t)$, if $G'$ is isomorphic to a graph in $\mathcal{E}_t$, then $G$ is also isomorphic to a graph in $\mathcal{E}_t$. In conclusion, $G$ is isomorphic to a graph in $\mathcal{E}_t$.

The proof of Theorem \ref{thm2} is complete. \qed

%%%%%%%%%%%%%%%%%%%%%%%%%%%%%%%%%%%%%%%

\section{Maximal independent sets in triangle-free graphs with a given matching number}\label{sec4}

Let 
\[A_n:=\begin{cases}
	\frac{n}{2}K_2\ &\text{for $n$ even},\\
	C_5\cup\frac{n-5}{2}K_2\ &\text{for $n$ odd}.
\end{cases}\]

\begin{letterthm}\cite{Hujtera1993}\label{thmD}
For any triangle-free graph $G$ of order $n\geq 4$, 
\[mis(G)\leq\begin{cases}
	2^{n/2}\ &\text{for\ } n \ \text{even},\\
	5\cdot 2^{(n-5)/2}\ &\text{for\ } n \ \text{odd},
\end{cases}\]
with equality if and only if $G\cong A_n$.
\end{letterthm}

Let
\begin{align*}
	\begin{split}
		\mathcal{M}_t:= \left \{
		\begin{array}{ll}
			\{\frac{t}{2}C_5\cup rK_1: r\in N\}                  & \text{for}\ t\ \text{even}, \\
			\{K_{1,\ell}\cup \frac{t-1}{2}C_5\cup sK_1:\ \ell\in Z^{+}\ \text{and}\ s\in N\} & \text{for}\ t\ \text{odd}.
		\end{array}
		\right.
	\end{split}
\end{align*}

\begin{thm}\label{thm3}
Let $G$ be a triangle-free graph with matching number $t$. Then
\[mis(G)\leq m(t):=\begin{cases}
	5^{t/2}\ &\text{for\ } t \ \text{even},\\
	2\cdot 5^{(t-1)/2}\ &\text{for\ } t \ \text{odd},
\end{cases}\]
with equality if and only if $G$ is isomorphic to a graph in $\mathcal{M}_t$.
\end{thm}

\pf Prove by induction on $t$.

When $t=1$, $G$ has only one non-trivial component, and this component can only be a star. Thus
\[mis(G)=2=m(1),\]
and $G$ is isomorphic to a graph in $\{K_{1,\ell}\cup sK_1:\ \ell\in Z^{+}\ \text{and}\ s\in N\}.$

Consider the case when $t=2$. If there is a vertex $v$ in $G$ that is saturated by all maximum matchings, we construct a new graph $G'$ from $G$ by adding a new vertex and joining the new vertex to $v$. By Lemmas \ref{lem3.2} and \ref{lem3.3}, $\mu(G'-v)=1$, $\mu(G'-N_{G'}[v])\leq 1$, and
\[mis(G)\leq mis(G')=mis(G'-v)+mis(G'-N_{G'}[v])\leq 2+2<5=m(2).\]

If there is no vertex in $G$ that is saturated by all maximum matchings, then by Corollary \ref{coro1}, each component of $G$ is factor-critical.
Since $G$ is triangle-free, $G$ has only one non-trivial component, and this component is isomorphic to $C_5$. 

In conclusion, when $t=2$, $mis(G)\leq m(2)$, with equality if and only if $G$ is isomorphic to a graph in $\{C_5\cup rK_1: r\in N\}$.

Assume that $t\geq 3$, and for any integer $t'<t$, if $H$ is a triangle-free graph with matching number $t'$, then $mis(H)\leq m(t')$, with equality if and only if $H$ is isomorphic to a graph in $\mathcal{M}_{t'}$.

Now, let $G_1$, $\cdots$, $G_k$ be non-trivial components of $G$, and let $\mu(G_1)=t_1$, $\cdots$, $\mu(G_k)=t_k$. Clearly, $t_1+\cdots+t_k=t$.

We will prove that if $k=1$, then 
\[mis(G)=mis(G_1)<m(t).\]

If there is a vertex $v$ in $G_1$ that is saturated by all maximum matchings, we construct a new graph $G'_1$ from $G_1$ by adding a new vertex and joining the new vertex to $v$. By Lemma \ref{lem3.2}, $\mu(G'_1-v)=t-1$ and $\mu(G'_1-N_{G'_1}[v])\leq t-1$. 

If $t$ is even, then
\begin{align*}
mis(G_1)&\leq mis(G'_1)\\
&=mis(G'_1-v)+mis(G'_1-N_{G'_1}[v])\\
&\leq 2\cdot 5^{(t-2)/2}+ 2\cdot 5^{(t-2)/2}\\
&<5^{t/2}=m(t).
\end{align*}

If $t$ is odd, then
\begin{align*}
	mis(G_1)&\leq mis(G'_1)\\
	&=mis(G'_1-v)+mis(G'_1-N_{G'_1}[v])\\
	&\leq 5^{(t-1)/2}+ 5^{(t-1)/2}\\
	&=2\cdot 5^{(t-1)/2}=m(t).
\end{align*}
Moreover, if $mis(G_1)=m(t)$, then by inductive hypothesis, $(G'_1-v)\cong \frac{t-1}{2}C_5\cup r_1K_1$ for some natural number $r_1$ and $(G'_1-N_{G'_1}[v])\cong \frac{t-1}{2}C_5\cup r_2K_1$ for some natural number $r_2$, it is impossible (since $G_1$ is connected). Thus, 
\[mis(G_1)<m(t).\]

If there is no vertex in $G_1$ that is saturated by all maximum matchings, then by Corollary \ref{coro1}, $G_1$ is factor-critical, which implies that $|V(G_1)|=2t+1\geq7$. By Theorem \ref{thmD},
\[mis(G_1)<5\cdot 2^{(t-2)}\leq 2\cdot 5^{(t-1)/2}\leq m(t).\] 

We have proved that if $k=1$, then $mis(G)=mis(G_1)<m(t).$ We have actually also proven that for any connected triangle-free graph $F$ with matching number $t\geq 3$, $mis(F)<m(t)$.

Now, suppose that $k\geq 2$. If there are at least two non-trivial components of $G$ with matching number $1$, then by inductive hypothesis,
\[mis(G)=2\cdot 2\cdot m(t-2)\leq 2\cdot 2\cdot 5^{(t-2)/2}<2\cdot 5^{(t-1)/2}\leq m(t).\]

Based on the previous analysis and conclusions, we have
\[mis(G)=\prod_{i=1}^k mis(G_i)\leq \prod_{i=1}^k m(t_i)\leq m(t).\]
Furthermore, if $mis(G)=m(t)$, then for each $i\in\{1,\cdots,k\}$, $t_i\leq 2$, and there is at most one non-trivial component of $G$ that is isomorphic to a star, and the other non-trivial components are isomorphic to $C_5$. In conclusion, if $mis(G)=m(t)$, then $G$ is isomorphic to a graph in $\mathcal{M}_t$. 

The proof of Theorem \ref{thm3} is complete. \qed

%%%%%%%%%%%%%%%%%%%%%%%%%%%%%%%%%%%%%%%

\section{Maximal independent sets in connected triangle-free graphs with a given matching number}\label{sec5}

When $n$ is even and $n\geq 10$, let 
\[B_n:=2C_5\cup\frac{n-10}{2}K_2.\]

\begin{figure}[H]
	\center
	\begin{tikzpicture}
		[line width = 1pt, scale=0.6,
		empty/.style = {circle, draw, fill = white, inner sep=0mm, minimum size=2mm}, full/.style = {circle, draw, fill = black, inner sep=0mm, minimum size=2mm}, full1/.style = {circle, draw, fill = black, inner sep=0mm, minimum size=0.5mm}]
		
		\node [empty] (a) at (1,1) {};
		\node [empty] (b) at (5,1) {};
		\node [empty] (c) at (-1,-1) {};
		\node [empty] (d) at (0.5,-1) {};
		\node [empty] (e) at (2,-1) {};
		\node [empty] (f) at (3.5,-1) {};
		\node [empty] (g) at (5,-1) {};

		\draw (a) -- (b);
		\draw (a) -- (c);
		\draw (a) -- (f);
		\draw (b) -- (g);
		\draw (c) -- (d);
		\draw (d) -- (e);
		\draw (e) -- (f);
				
		\path (2,-3.5) node(text1)[right]{$H_7$};

		%%%%%%%%%%%%%%%%%%%%%%%%%%%%%%%%%%%%%%%%%%
		
		\node [empty] (a) at (13,2) {};
		\node [empty] (b) at (11,0.5) {};
		\node [empty] (c) at (11,-1) {};
		\node [empty] (d) at (15,0.5) {};
		\node [empty] (e) at (15,-1) {};

		\node [label=above:$\cdots$] at (13,-1) {};

		\draw (a) -- (b);
		\draw (b) -- (c);
		\draw (a) -- (d);
		\draw (d) -- (e);
		
		\draw[decorate,decoration={brace,mirror,amplitude=10pt},xshift=0pt,yshift=-4pt]
		(11,-1.2) -- (15,-1.2) node [black,midway,yshift=-0.7cm] {$r$};
			
		\path (12,-3.5) node(text1)[right]{$T_{2r+1}$};
		
	\end{tikzpicture}
	\caption{Graphs $H_7$ and $T_{2r+1}$.}
	\label{fig3}
\end{figure}
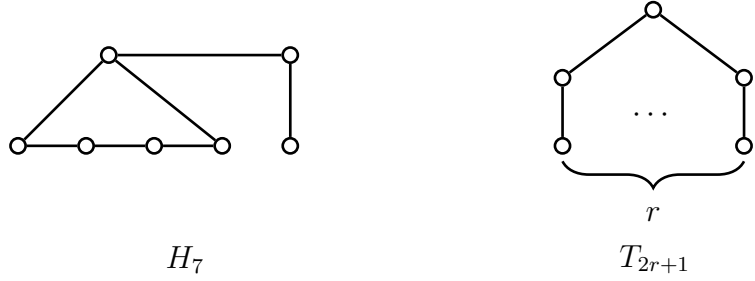

When $n$ is odd, let 
\[D_n:=H_7\cup\frac{n-7}{2}K_2\ \text{or}\ T_{2r+1}\cup \frac{n-2r-1}{2}K_2 \text{\ with } 0\leq r\leq \frac{n-1}{2},\]
where $H_7$ and $T_{2r+1}$ are shown in Figure \ref{fig3}.

\begin{letterthm}\cite{Chang1999}\label{thmE}
	Let $G$ be a triangle-free graph $G$ of order $n$. If $G\not\cong A_n$ and $G\not\cong B_n$, then
\[mis(G)\leq q(n):=\begin{cases}
			3\cdot 2^{\frac{n-4}{2}}          & \text{for $n$ even}, \\
			2^\frac{n-1}{2}             & \text{for $n$ odd}.
\end{cases}\]
Furthermore, when $n$ is odd, $mis(G)=q(n)$ if and only if $G\cong D_n$.
\end{letterthm}

\begin{figure}[H]
	\center
	\begin{tikzpicture}
		[line width = 1pt, scale=0.6,
		empty/.style = {circle, draw, fill = white, inner sep=0mm, minimum size=2mm}, full/.style = {circle, draw, fill = black, inner sep=0mm, minimum size=2mm}, full1/.style = {circle, draw, fill = black, inner sep=0mm, minimum size=0.5mm}]
		\node [empty] (a) at (1,0) {};
		\node [empty,label = below:$u$] (b) at (3,0) {};
		\node [empty] (c) at (5,0) {};
		\node [empty,label = below:$v$] (d) at (7,0) {};
		\node [empty] (e) at (9,0) {};
		%\node [empty] (f) at (1,1) {};
		
		\node [label=above:$\vdots$] at (1,0.2) {};
		
			\draw[decorate,decoration={brace,amplitude=10pt},xshift=-4pt,yshift=0pt]
		(0.8,0) -- (0.8,2) node [midway,xshift=-0.6cm] {$\ell_1$};

		\node [empty] (g) at (1,2) {};

		\node [empty] (h) at (5,1.5) {};
		\node [label=above:$\vdots$] at (5,0) {};
		
			\draw[decorate,decoration={brace,mirror,amplitude=10pt},xshift=-4pt,yshift=0pt]
		(5.4,0) -- (5.4,1.5) node [midway,xshift=0.6cm] {$\ell_3$};

		\node [empty] (j) at (9,2) {};
		
		\node [label=above:$\vdots$] at (9,0.2) {};
		
					\draw[decorate,decoration={brace,mirror,amplitude=10pt},xshift=-4pt,yshift=0pt]
		(9.4,0) -- (9.4,2) node [midway,xshift=0.6cm] {$\ell_2$};

		\node [empty] (k) at (1.5,-2) {};
		\node [empty] (l) at (0.5,-3) {};
		\node [empty] (m) at (1,-4) {};
		\node [empty] (n) at (2,-4) {};
		\node [empty] (o) at (2.5,-3) {};

		\node [empty] (p) at (4.5,-2) {};
		\node [empty] (q) at (3.5,-3) {};
		\node [empty] (r) at (4,-4) {};
		\node [empty] (s) at (5,-4) {};
		\node [empty] (t) at (5.5,-3) {};
		
		\node [label=above:$\cdots$] at (6.6,-3.8) {};
		
		\node [empty] (u) at (8.5,-2) {};
		\node [empty] (v) at (7.5,-3) {};
		\node [empty] (w) at (8,-4) {};
		\node [empty] (x) at (9,-4) {};
		\node [empty] (y) at (9.5,-3) {};
		%\node [full1] at (2.6,1) {};
		\draw (a) -- (b);
		\draw (b) -- (c);
		\draw (d) -- (c);
		\draw (d) -- (e);
		%\draw (f) -- (b);
		\draw (g) -- (b);
		\draw (h) -- (b);
		\draw (h) -- (d);
		\draw (j) -- (d);
		\draw (k) -- (l);
		\draw (m) -- (l);
		\draw (m) -- (n);
		\draw (o) -- (n);
		\draw (o) -- (k);
		\draw (b) -- (k);
		\draw (p) -- (q);
		\draw (r) -- (q);
		\draw (r) -- (s);
		\draw (t) -- (s);
		\draw (t) -- (p);
		\draw (b) -- (p);
		\draw (u) -- (v);
		\draw (w) -- (v);
		\draw (w) -- (x);
		\draw (y) -- (x);
		\draw (y) -- (u);
		\draw (b) -- (u);
		
		\draw[decorate,decoration={brace,mirror,amplitude=10pt},xshift=0pt,yshift=-4pt]
		(1,-4.5) -- (8.5,-4.5) node [black,midway,yshift=-0.7cm] {$\frac{t-2}{2}$};

		\path (3.2,-7) node(text1)[right]{For $t$ even};

		%%%%%%%%%%%%%%%%%%%%%%%%%%%%%%%%%%%%%%%%%%
		
		\node [empty] (a) at (13,0) {};
		\node [empty,label=right:$w$] (b) at (15,0) {};
		%\node [empty] (f) at (13,1) {};
		\node [label=above:$\vdots$] at (13,0.2) {};
		
		\node [empty] (g) at (13,2) {};

		\node [empty] (k) at (13.5,-2) {};
		\node [empty] (l) at (12.5,-3) {};
		\node [empty] (m) at (13,-4) {};
		\node [empty] (n) at (14,-4) {};
		\node [empty] (o) at (14.5,-3) {};

		\node [empty] (p) at (16.5,-2) {};
		\node [empty] (q) at (15.5,-3) {};
		\node [empty] (r) at (16,-4) {};
		\node [empty] (s) at (17,-4) {};
		\node [empty] (t) at (17.5,-3) {};
		
		\node [label=above:$\cdots$] at (18.6,-3.8) {};

		\node [empty] (u) at (20.5,-2) {};
		\node [empty] (v) at (19.5,-3) {};
		\node [empty] (w) at (20,-4) {};
		\node [empty] (x) at (21,-4) {};
		\node [empty] (y) at (21.5,-3) {};
		%\node [full1] at (2.6,1) {};
		\draw (a) -- (b);
	%	\draw (b) -- (f);
		\draw (b) -- (g);
		\draw (k) -- (l);
		\draw (m) -- (l);
		\draw (m) -- (n);
		\draw (o) -- (n);
		\draw (o) -- (k);
		\draw (b) -- (k);
		\draw (p) -- (q);
		\draw (r) -- (q);
		\draw (r) -- (s);
		\draw (t) -- (s);
		\draw (t) -- (p);
		\draw (b) -- (p);
		\draw (u) -- (v);
		\draw (w) -- (v);
		\draw (w) -- (x);
		\draw (y) -- (x);
		\draw (y) -- (u);
		\draw (b) -- (u);
		
		\draw[decorate,decoration={brace,mirror,amplitude=10pt},xshift=0pt,yshift=-4pt]
		(13,-4.5) -- (20.5,-4.5) node [black,midway,yshift=-0.7cm] {$\frac{t-1}{2}$};
		
		\path (15.6,-7) node(text1)[right]{For $t$ odd};

	\end{tikzpicture}
	\caption{Two examples of graphs in $\mathcal{G}_t$}
	\label{fig4}
\end{figure}
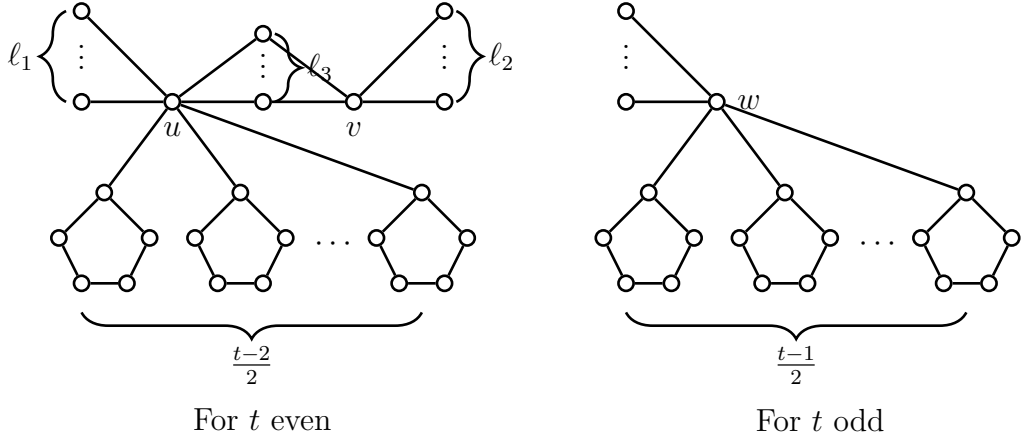

For a given positive integer $t$, we introduce a set of graphs $\mathcal{G}_t$. When $t$ is even, the graphs in $\mathcal{G}_t$ can be constructed through the following steps:

\begin{itemize}

\item[(1)] Construct three sets of vertices $L_1$, $L_2$, and $L_3$.

\item[(2)] Add two vertices $u$ and $v$, and connect each vertex in $L_1$ to $u$, each vertex in $L_2$ to $v$, and each vertex in $L_3$ to both $u$ and $v$.

 (The class of graphs obtained at this step is referred to as $\mathcal{P}$.)

\item[(3)] Add $\frac{t-2}{2}$ vertex-disjoint $C_5$s, and connect exactly one vertex of each $C_5$ to $u$.
\end{itemize}
\noindent When $t$ is odd, the graphs in $\mathcal{G}_t$ can be constructed through the following steps:
\begin{itemize}
\item[(1)] Form a star $K_{1,r}$, where $r \geq 1$, and denote the central vertex of the star as $w$.

\item[(2)] Add $\frac{t-1}{2}$ vertex-disjoint $C_5$s, and connect exactly one vertex of each $C_5$ to $w$.
\end{itemize}
Two examples of graphs in $\mathcal{G}_t$ are shown in Figure \ref{fig4}.

\begin{figure}[H]
	\center
	\begin{tikzpicture}
		[line width = 1pt, scale=0.7,
		empty/.style = {circle, draw, fill = white, inner sep=0mm, minimum size=2mm}, full/.style = {circle, draw, fill = black, inner sep=0mm, minimum size=2mm}, full1/.style = {circle, draw, fill = black, inner sep=0mm, minimum size=0.5mm}]
		\node [empty] (a) at (1,1) {};
		
		\node [empty,label = below:$u$] (b) at (3,0) {};
		\node [empty] (c) at (5,0) {};
		\node [empty,label = below:$v$] (d) at (7,0) {};
		\node [empty] (e) at (9,0) {};

		%\node [empty] (f) at (1,1) {};
		%\node [empty] (g) at (1,2) {};
		\node [empty] (h) at (1,0) {};
		\node [empty] (i) at (1,-1) {};

		\node [empty] (j) at (5,1) {};
		\node [empty] (k) at (5,-1) {};
		%\node [full1] at (2.6,1) {};
		\node [empty] (l) at (9,1) {};
		\node [empty,label=right:$w$] (m) at (3,3) {};
		
		\node [empty] (n) at (2,4) {};

		%\node [empty] (o) at (3,4) {};
		\node [label=above:$\cdots$] at (3,3.4) {};

		\node [empty] (p) at (4,4) {};
		\draw (a) -- (b);
		\draw (b) -- (c);
		\draw (d) -- (c);
		\draw (d) -- (e);
		%\draw (f) -- (b);
	%	\draw (g) -- (b);
		\draw (h) -- (b);
		\draw (i) -- (b);
		\draw (b) -- (j);
		\draw (j) -- (d);
		\draw (k) -- (b);
		\draw (k) -- (d);
		\draw (d) -- (l);
		\draw (n) -- (m);
	%	\draw (o) -- (m);
		\draw (p) -- (m);
		\draw (a) -- (m);
	%	\draw (f) -- (m);
	%	\draw (g) -- (m);
	%	\draw (h) -- (m);
		\draw (c) -- (m);
		\draw (j) -- (m);
		\draw (l) -- (m);
		\draw (k) -- (m);
	%	\path (3,-3) node(text1)[right]{A graph in $\mathcal{Q}_3$};
		
%%%%%%%%%%%%%%%%%%%%%%%%%%%%%%%%%%%%%%%

	\node [empty] (a) at (12,1) {};

\node [empty] (b) at (14,0) {};
%\node [empty] (c) at (16,0) {};
\node [empty] (d) at (16,0) {};
\node [empty] (e) at (18,1) {};
\node [empty] (e1) at (18,0) {};
\node [empty] (e2) at (18,-1) {};

%\node [empty] (f) at (1,1) {};
%\node [empty] (g) at (1,2) {};
\node [empty] (h) at (12,0) {};
\node [empty] (i) at (12,-1) {};

%\node [empty] (j) at (16,1) {};
%\node [empty] (k) at (16,-1) {};
%\node [full1] at (2.6,1) {};
\node [empty] (l) at (18,2) {};
\node [empty,label=right:$w$] (m) at (14,3) {};

\node [empty] (n) at (13,4) {};

%\node [empty] (o) at (3,4) {};
\node [label=above:$\cdots$] at (14,3.4) {};

\node [empty] (p) at (15,4) {};
\draw (a) -- (b);
%\draw (b) -- (c);
%\draw (d) -- (c);
\draw (d) -- (e);
\draw (d) -- (e1);
\draw (d) -- (e2);
%\draw (f) -- (b);
%	\draw (g) -- (b);
\draw (h) -- (b);
\draw (i) -- (b);
%\draw (b) -- (j);
%\draw (j) -- (d);
%\draw (k) -- (b);
%\draw (k) -- (d);
\draw (d) -- (l);
\draw (n) -- (m);
%	\draw (o) -- (m);
\draw (p) -- (m);
\draw (a) -- (m);
%	\draw (f) -- (m);
%	\draw (g) -- (m);
%	\draw (h) -- (m);
%\draw (c) -- (m);
%\draw (j) -- (m);
\draw (l) -- (m);
%\draw (k) -- (m);
\draw (h)--(m);
	\end{tikzpicture}
	\caption{Two examples of graphs in $\mathcal{Q}_3$}
	\label{fig5}
\end{figure}
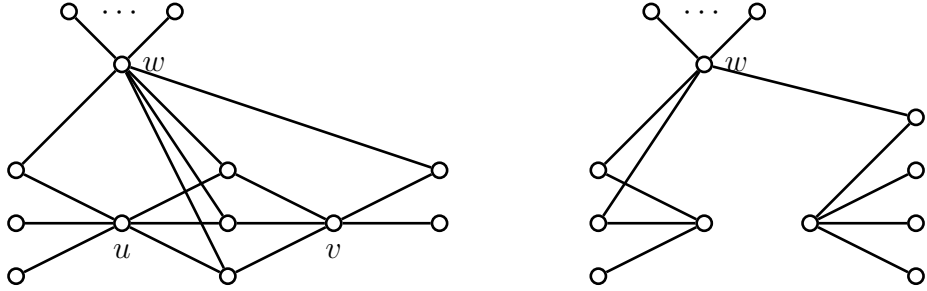

Next, we introduce a set of graphs $\mathcal{Q}_3$. The graphs in $\mathcal{Q}_3$ can be constructed through the following steps:

\begin{itemize}
	
   \item[(1)] Choose a graph $H$ from $\mathcal{P}\cup \{K_{1,r_1}\cup K_{1,r_2}:\ r_1\geq 2,\ r_2\geq 2\}$.
   
   	\item[(2)] Construct a star $K_{1,s}$, where $s\geq 1$, and denote the central vertex of the star as $w$.
	
	\item[(3)] Connect $w$ to some neighbors of the supported vertices in $H$, ensuring that these supported vertices are still supported vertices in the new graph.
\end{itemize}
Two examples of graphs in $\mathcal{Q}_3$ are shown in Figure \ref{fig5}.

Finally, we introduce a set of graphs $\mathcal{Q}_4$. The graphs in $\mathcal{Q}_4$ can be constructed through the following steps:

\begin{itemize}
	
	\item[(1)] Choose a graph $H$ from 
\[\mathcal{Q}_3\cup\{F\cup K_{1,r}:\ F\in \mathcal{P},\ r\geq 2\}\cup \{K_{1,r_1}\cup K_{1,r_2}\cup K_{1,r_3}:\ r_1\geq 2,\ r_2\geq 2,\ r_3\geq2\}.\]
	
	\item[(2)] Construct a star $K_{1,q}$, where $q\geq 1$, and denote the central vertex of the star as $w$,
	
	\item[(3)] Connect $w$ to some neighbors of the supported vertices in $H$, ensuring that these supported vertices are still supported vertices in the new graph.
\end{itemize}

Let
\begin{align*}
	\begin{split}
		\mathcal{F}_t:= \left \{
		\begin{array}{ll}
			\{K_{1,r}:r\in Z^+\}                  &\text{for $t=1$}, \\
			\{C_5\}                  &\text{for $t=2$}, \\
			\mathcal{G}_3\cup\mathcal{Q}_3                &\text{for $t=3$}, \\
			\mathcal{G}_4\cup\mathcal{Q}_4                &\text{for $t=4$}, \\
			\mathcal{G}_t  &\text{for $t\geq5$}.
		\end{array}
		\right.
	\end{split}
\end{align*}

\begin{thm}\label{thm4}
For any connected triangle-free graph $G$ with matching number $t$, 
\[mis(G)\leq f(t):=\begin{cases}
	5           & \text{for $t=2$},\\
	2\cdot(5^{\frac{t-2}{2}}+3^{\frac{t-2}{2}}) & \text{for $t\geq 4$ even}, \\
	5^{\frac{t-1}{2}}+3^{\frac{t-1}{2}}    & \text{for $t$ odd},
\end{cases}\]
with equality if and only if $G$ is isomorphic to a graph in $\mathcal{F}_t$.
\end{thm}

\begin{rem}
	It is not hard to check that $f(t)$ is a strictly increasing function with respect to positive integer $t$.
\end{rem}

\pf Let $\phi(t)$ be the maximum number of MISs among all connected triangle-free graphs with matching number $t$.

For any positive integer $t$ and any graph $G\in \mathcal{F}_t$, by direct calculations, 
\[mis(G)=f(t),\]
which implies that $\phi(t)\geq f(t)$.

We prove by induction on $t$ that $\phi(t)\leq f(t)$, and any connected triangle-free graph with matching number $t$ and $f(t)$ MISs is isomorphic to a graph in $\mathcal{F}_t$.

When $t=1$, it is easy to see that any connected triangle-free graph $G$ with matching number $1$ is a star and has $f(1)=2$ MISs.

Suppose that $t\geq2$, and for any positive integer $t'<t$, $\phi(t')\leq f(t')$, and any connected triangle-free graph with matching number $t'$ and $f(t')$ MISs is isomorphic to a graph in $\mathcal{F}_{t'}$.

Let $G$ be a connected triangle-free graph with matching number $t$ and $\phi(t)$ MISs. We distinguish the following cases.

\noindent {\bf Case 1}. There is no vertex in $G$ that is saturated by all maximum matchings.

In this case, by Corollary \ref{coro1} and the fact that $G$ is connected, $G$ is factor-critical, which implies that $|V(G)|=2t+1$.

If $t=2$, then $G\cong C_5$ and $mis(G)=f(2)=5$.

If $t\geq 3$, then $G\not\cong A_{2t+1}$, $G\not\cong B_{2t+1}$ and $G\not\cong D_{2t+1}$. By Theorem \ref{thmE}, 
\[mis(G)< q(2t+1)=2^{t}\leq f(t).\]

\noindent {\bf Case 2.} There exists a vertex $v$ in $G$ that is saturated by all maximum matchings.

We obtain a graph $G'$ from $G$ by adding a new vertex and joining the new vertex to $v$. By Lemmas \ref{lem3.1} and \ref{lem3.3}, $\mu(G')=\mu(G)=t$ and $mis(G')\geq mis(G)=\phi(t)$, which implies that $mis(G')=\phi(t)$. Furthermore,
\[mis(G')=mis(G'-v)+mis(G'-N_{G'}[v]),\]
and
\[\mu(G'-v)=t-1.\]

When $t=2$, $mis(G'-N_{G'}[v])\leq mis(G'-v)\leq 2$, which implies that 
\[mis(G')\leq 4< f(2).\]
Furthermore, $mis(G')=4$ if and only if $G'$ is isomorphic to a graph in $\mathcal{P}$.

%%%%%%%%%%%%%%%%%%%
When $t=3$, by Theorem \ref{thm3}, 
\[mis(G'-N_{G'}[v])\leq mis(G'-v)\leq m(2)=5.\]
If $mis(G'-v)=5$, then $G'-v\cong C_5\cup rK_1$ for some $r\geq 1$, and 
$mis(G'-N_{G'}[v])\leq 3$, with equality if and only if $G'-N_{G'}[v]\cong P_4$. It implies that if $mis(G'-v)=5$, then $mis(G')\leq 8$, with equality if and only if $G'$ is isomorphic to a graph in $\mathcal{G}_3$. 

If $mis(G'-N_{G'}[v])=4$ and $mis(G'-v)=4$, then $G'-v$ is isomorphic to a graph in $\mathcal{P}\cup \{K_{1,r_1}\cup K_{1,r_2}:\ r_1\geq 2,\ r_2\geq 2\}$, and $G'$ is isomorphic to a graph in $\mathcal{Q}_3$.

In conclusion, $mis(G')\leq f(3)=8$, with equality if and only if $G'$ is isomorphic to a graph in $\mathcal{F}_3=\mathcal{G}_3\cup \mathcal{Q}_3$.

%%%%%%%%%%%%%%%%%%%%%%%%%%%%%%

When $t=4$, by Theorem \ref{thm3}, 
\[mis(G'-N_{G'}[v])\leq mis(G'-v)\leq m(3)=10.\]
If $mis(G'-v)=10$, then $G'-v\cong K_{1,\ell}\cup C_5\cup sK_1$ for some $\ell\geq 1$ and $s\geq 1$, and $mis(G'-N_{G'}[v])\leq 6$, with equality if and only if $G'-N_{G'}[v]\cong P_4\cup K_{1,\ell'}$ for some $1\leq \ell'<\ell$. It implies that if $mis(G'-v)=10$, then $mis(G')\leq 16$, with equality if and only if $G'$ is isomorphic to a graph in $\mathcal{G}_4$. 

If $mis(G'-v) = 9$, then by inductive hypothesis, $G'-v$ is disconnected and contains two non-trivial components, each of which has $3$ MISs. However, since $\mu(G'-v) = 3$, this scenario is impossible.

If $mis(G'-v)=8$ and $mis(G'-N_{G'}[v])=8$, then $G'-v$ is isomorphic to a graph in 
\[\mathcal{Q}_3\cup\{F\cup K_{1,r}:\ F\in \mathcal{P},\ r\geq 2\}\cup \{K_{1,r_1}\cup K_{1,r_2}\cup K_{1,r_3}:\ r_1\geq 2,\ r_2\geq 2,\ r_3\geq2\},\]
and $G'$ is isomorphic to a graph in $\mathcal{Q}_4$.

In conclusion, $mis(G')\leq f(4)=16$, with equality if and only if $G'$ is isomorphic to a graph in $\mathcal{F}_4=\mathcal{G}_4\cup \mathcal{Q}_4$.\\

%%%%%%%%%%%%%%%%%%%%%%%%%%%%%%%

Now, we suppose that $t\geq 5$.

\begin{cl}\label{cl4.1}
In $G'-v$, there is no component $H$ such that $\mu(H)$ is even and $\mu(H)\geq 4$.
\end{cl}

\noindent{\it Proof of Claim \ref{cl4.1}.} 
Suppose, to the contrary, that $H$ is a component of $G'-v$ such that $\mu(H)$ is even and $\mu(H)=k\geq 4$. Let $\alpha:=mis(H)$ and $\beta:=mis(H-(N_{G'}(v)\cap V(H)))$. Let $G_1:=G'-v-V(H)$ and $G_2:=G'-\left[N_{G'}[v]\cup V(H)\right]$. Clearly, $V(G_1)\geq 1$. Thus, 
\[mis(G')=mis(G'-v)+mis(G'-N_{G'}[v])=\alpha\cdot mis(G_1)+\beta\cdot mis(G_2).\]
By inductive hypothesis, 
\[\beta\leq\alpha\leq 2\cdot(5^{\frac{k-2}{2}}+3^{\frac{k-2}{2}}).\]

We obtain a new graph $G''$ from $G'$ by first removing all vertices in $H$ from $G'$, then adding $\frac{k}{2}$ new $C_5$s, and finally connecting the vertex $v$ to exactly one vertex in each of these newly added $C_5$s. It is easy to see that $G''$ is a connected triangle-free graph with matching number $t$. Since $mis(G_1)\geq mis(G_2)\geq 1$ and $5^\frac{k}{2}+3^\frac{k}{2}>4\cdot(5^{\frac{k-2}{2}}+3^{\frac{k-2}{2}})$, we have
\begin{align*}
	\phi(t)=&mis(G')=\alpha\cdot mis(G_1)+\beta\cdot mis(G_2)\\
	\leq & 2\cdot(5^{\frac{k-2}{2}}+3^{\frac{k-2}{2}})\cdot mis(G_1)+2\cdot(5^{\frac{k-2}{2}}+3^{\frac{k-2}{2}})\cdot mis(G_2)\\
	<& 2\cdot(5^{\frac{k-2}{2}}+3^{\frac{k-2}{2}})\cdot mis(G_1)+[5^\frac{k}{2}-2\cdot(5^{\frac{k-2}{2}}+3^{\frac{k-2}{2}})]\cdot mis(G_2)+3^{\frac{k}{2}}\cdot mis(G_2)\\
	\leq & 2\cdot(5^{\frac{k-2}{2}}+3^{\frac{k-2}{2}})\cdot mis(G_1)+[5^\frac{k}{2}-2\cdot(5^{\frac{k-2}{2}}+3^{\frac{k-2}{2}})]\cdot mis(G_1)+3^{\frac{k}{2}}\cdot mis(G_2)\\
	= & 5^\frac{k}{2}\cdot mis(G_1)+3^\frac{k}{2}\cdot mis(G_2)\\
	=&mis(G''),
\end{align*}
which contradicts that $\phi(t)$ is the maximum number of MISs among all connected triangle-free graphs with matching number $t$. \qed

\begin{cl}\label{cl4.2}
	If there exists at least two components with matching number $2$ in $G'-v$, then each of these components is isomorphic to $C_5$, and there is exactly one vertex in each that is adjacent to $v$.
\end{cl}

\noindent{\it Proof of Claim \ref{cl4.2}.} Let $H_1$ and $H_2$ be two components with matching number $2$ in $G'-v$. Let $\alpha_1:=mis(H_1)$, $\beta_1:=mis(H_1-(N_{G'}(v)\cap V(H_1)))$, $\alpha_2:=mis(H_2)$, and $\beta_2:=mis(H_2-(N_{G'}(v)\cap V(H_2)))$. Let $G_3:=G'-v-V(H_1)-v(H_2)$ and $G_4:=G'-\left[N_{G'}[v]\cup V(H_1)\cup V(H_2)\right]$. Thus,
\[mis(G')=\alpha_1\alpha_2\cdot mis(G_3)+\beta_1\beta_2\cdot mis(G_4).\]
Based on the previous analysis for the case when $t=2$, for every $1 \leq i \leq 2$, the best possible values for the pair $(\alpha_i, \beta_i)$ are $(5,3)$ and $(4,4)$. Therefore, $mis(G')$ is maximized only when both $(\alpha_1, \beta_1)$ and $(\alpha_2, \beta_2)$ are $(5,3)$. Thus, both $H_1$ and $H_2$ are isomorphic to $C_5$, and in both $H_1$ and $H_2$, there is only one vertex adjacent to $v$. \qed

\begin{cl}\label{cl4.3}
If there exists a component $H_1$ in $G'-v$ such that $\mu(H_1)$ is odd and $\mu(H_1)\geq 3$, then $H_1$ is the only component in $G'-v$ whose matching number is odd. 
\end{cl}

\noindent{\it Proof of Claim \ref{cl4.3}.} Suppose, to the contrary, that $H_2$ is another component in $G'-v$ such that $\mu(H_2)$ is odd. Let $\mu(H_1)=t_1$ and $\mu(H_2)=t_2$. 

Let $\alpha_1:=mis(H_1)$, $\beta_1:=mis(H_1-(N_{G'}(v)\cap V(H_1)))$, $\alpha_2:=mis(H_2)$, $\beta_2:=mis(H_2-(N_{G'}(v)\cap V(H_2)))$. Let $G_3:=G'-v-V(H_1)-v(H_2)$ and $G_4:=G'-\left[N_{G'}[v]\cup V(H_1)\cup V(H_2)\right]$. Clearly, $mis(G_3)\geq mis(G_4)\geq1$.

We obtain a new graph $G''$ from $G'$ by first removing all vertices in $H_1$ and $H_2$ from $G'$, then adding $\frac{t_1+t_2}{2}$ new $C_5$s, and finally connecting the vertex $v$ to exactly one vertex in each of these newly added $C_5$s. It is easy to see that $G''$ is a connected triangle-free graph with matching number $t$. 

Thus, by inductive hypothesis and the following inequality
\[5^{\frac{t_1+t_2}{2}}+3^{\frac{t_1+t_2}{2}}>2\cdot (5^{\frac{t_1-1}{2}}+3^{\frac{t_1-1}{2}})(5^{\frac{t_2-1}{2}}+3^{\frac{t_2-1}{2}}),\]
we have
\begin{align*}
	\phi(t)=&mis(G')=mis(G'-v)+mis(G'-N_{G'}[v])\\
	=&\alpha_1\alpha_2\cdot mis(G_3)+\beta_1\beta_2\cdot mis(G_4)\\
	\leq & (5^{\frac{t_1-1}{2}}+3^{\frac{t_1-1}{2}})(5^{\frac{t_2-1}{2}}+3^{\frac{t_2-1}{2}})\cdot mis(G_3)+(5^{\frac{t_1-1}{2}}+3^{\frac{t_1-1}{2}})(5^{\frac{t_2-1}{2}}+3^{\frac{t_2-1}{2}})\cdot mis(G_4)\\
	<&(5^{\frac{t_1-1}{2}}+3^{\frac{t_1-1}{2}})(5^{\frac{t_2-1}{2}}+3^{\frac{t_2-1}{2}})\cdot mis(G_3)\\
	&+[5^{\frac{t_1+t_2}{2}}-(5^{\frac{t_1-1}{2}}+3^{\frac{t_1-1}{2}})(5^{\frac{t_2-1}{2}}+3^{\frac{t_2-1}{2}})]\cdot mis(G_4)+3^{\frac{t_1+t_2}{2}} \cdot mis(G_4)\\
	\leq&(5^{\frac{t_1-1}{2}}+3^{\frac{t_1-1}{2}})(5^{\frac{t_2-1}{2}}+3^{\frac{t_2-1}{2}})\cdot mis(G_3)\\
	&+[5^{\frac{t_1+t_2}{2}}-(5^{\frac{t_1-1}{2}}+3^{\frac{t_1-1}{2}})(5^{\frac{t_2-1}{2}}+3^{\frac{t_2-1}{2}})]\cdot mis(G_3)+3^{\frac{t_1+t_2}{2}} \cdot mis(G_4)\\
	=& 5^\frac{t_1+t_2}{2}mis(G_3)+ 3^\frac{t_1+t_2}{2}mis(G_4)\\
	=& mis(G'')
\end{align*}
which contradicts that $\phi(t)$ is the maximum number of MISs among all connected triangle-free graphs with matching number $t$. \qed\\

Let $R_1,\cdots,R_k$ be non-trivial components in $G'-v$. For each $1\leq i\leq k$, let $\mu(R_i)=r_i$. Without loss of generality, suppose $r_1\geq\cdots\geq r_k$.

\begin{cl}\label{cl4.4}
	$r_1\geq 2$.
\end{cl}
\noindent{\it Proof of Claim \ref{cl4.4}.} Suppose, to the contrary, that $r_1=1$. Since $t\geq 5$, there is at least four components with matching number $1$ in $G'-v$. Let $H=R_1\cup R_2\cup R_3\cup R_4$.

We construct a new graph $G''$: firstly remove all vertices in $H$ from $G'$, subsequently add two new $C_5$s, and finally connect the vertex $v$ to exactly one vertex in each of the two newly added $C_5$s. It is easy to see that $G''$ is a connected triangle-free graph with matching number $t$. 

Let $\alpha:=mis(H)$ and $\beta:=mis(H-(N_{G'}(v)\cap V(H)))$. Let $G_1:=G'-v-V(H)$ and $G_2:=G'-\left[N_{G'}[v]\cup V(H)\right]$. We have 
\begin{align*}
	\phi(t)=&mis(G')=mis(G'-v)+mis(G'-N_{G'}[v])\\
	=&\alpha\cdot mis(G_1)+\beta\cdot mis(G_2)\\
	\leq & 2^4\cdot mis(G_1)+2^4\cdot mis(G_2)\\
	<& 5^2\cdot mis(G_1)+ 3^2\cdot mis(G_2)\ \ \text{(since $mis(G_1)\geq mis(G_2))$}\\
	=& mis(G''),
\end{align*}
which contradicts that $\phi(t)$ is the maximum number of MISs among all connected triangle-free graphs with matching number $t$. \qed

\begin{cl}\label{cl4.5}
If $r_1\geq 3$, then $G'-v\cong R_1\cup sK_1$ for some $s\geq 1$, moreover, $G'$ is isomorphic to a graph in $\mathcal{F}_t$.
\end{cl}

\noindent{\it Proof of Claim \ref{cl4.5}.} By Claims \ref{cl4.1} and \ref{cl4.3}, $r_1$ is odd, and there is no other component than $R_1$ in $G'-v$ whose matching number is odd. 

Suppose, to the contrary, that $r_2=2$. We construct a new graph $G''$: firstly remove all vertices in $R_1$ from $G'$, subsequently add $\frac{r_1-1}{2}$ new $C_5$s and a star $K_{1,\ell}$ for some $\ell\geq 2$, and finally connect the vertex $v$ to exactly one vertex in each of these newly added $C_5$s, and to a leaf of the star. It is easy to see that $G''$ is a connected triangle-free graph with matching number $t$. 

Let $\alpha:=mis(R_1)$ and $\beta:=mis(R_1-(N_{G'}(v)\cap V(R_1)))$. Let $G_1:=G'-v-V(R_1)$ and $G_2:=G'-\left[N_{G'}[v]\cup V(R_1)\right]$. 

By inductive hypothesis,
\begin{align*}
	\phi(t)=&mis(G')=mis(G'-v)+mis(G'-N_{G'}[v])\\
	=&\alpha\cdot mis(G_1)+\beta\cdot mis(G_2)\\
	\leq & (5^\frac{r_1-1}{2}+3^\frac{r_1-1}{2})\cdot mis(G_1)+(5^\frac{r_1-1}{2}+3^\frac{r_1-1}{2})\cdot mis(G_2)\\
	\leq & 2\cdot5^\frac{r_1-1}{2}\cdot mis(G_1)+2\cdot3^\frac{r_1-1}{2}\cdot mis(G_2)\ \ \text{(since $mis(G_1)\geq mis(G_2))$}\\
	=& mis(G'').
\end{align*}
Due to the definition of $\phi(t)$, we have $mis(G'')=\phi(t)$, which implies that $mis(G_1)=mis(G_2)$. But, if $mis(G'')=\phi(t)$, then by Claim $\ref{cl4.2}$, $R_2\cong C_5$, which implies that it is impossible for $mis(G_1)=mis(G_2)$. This is a contradiction. Thus, $G'-v\cong R_1\cup sK_1$ for some $s\geq 1$.

Now, $t=r_1+1$ is even and $t\geq 6$. By inductive hypothesis,
\begin{align*}
	\phi(t)=mis(G')=&mis(G'-v)+mis(G'-N_{G'}[v])\\
	\leq&f(t-1)+f(t-1)=(5^{\frac{t-2}{2}}+3^{\frac{t-2}{2}})+(5^{\frac{t-2}{2}}+3^{\frac{t-2}{2}})\\
	\leq& f(t).
\end{align*}
Thus, $\phi(t)=f(t)$. We have $mis(G'-v)=mis(G'-N_{G'}[v])=f(t-1)$, and both $G'-v$ and $G'-N_{G'}[v]$ are isomorphic to graphs in $\mathcal{F}_{t-1}$.
It is not hard to check that $G'$ is isomorphic to a graph in $\mathcal{F}_t$. \qed

\begin{cl}\label{cl4.6}
	If $r_1=2$, then there exists at most one component with matching number $1$ in $G'-v$, moreover, $G'$ is isomorphic to a graph in $\mathcal{F}_t$.
\end{cl}

\noindent{\it Proof of Claim \ref{cl4.6}.} Without loss of generality, we suppose, for the sake of contradiction, that both $R_2$ and $R_3$ have matching number $1$.

We construct a new graph $G''$: firstly remove all vertices in $R_2$ and $R_3$ from $G'$, subsequently add a new $C_5$, and finally connect the vertex $v$ to exactly one vertex of the newly added $C_5$. Clearly, $G''$ is a connected triangle-free graph with matching number $t$. 

Let $H:=R_2\cup R_3$. Let $\alpha:=mis(H)$ and $\beta:=mis(H-(N_{G'}(v)\cap V(H)))$. Let $G_1:=G'-v-V(H)$ and $G_2:=G'-\left[N_{G'}[v]\cup V(H)\right]$. By inductive hypothesis,
\begin{align*}
	\phi(t)=&mis(G')=mis(G'-v)+mis(G'-N_{G'}[v])\\
	=&\alpha\cdot mis(G_1)+\beta\cdot mis(G_2)\\
	\leq & 4\cdot mis(G_1)+4\cdot mis(G_2)\\
	\leq & 5\cdot mis(G_1)+3\cdot mis(G_2)\ \ \text{(since $mis(G_1)\geq mis(G_2))$}\\
	=& mis(G'').
\end{align*}
Due to the definition of $\phi(t)$, we have $mis(G'')=\phi(t)$, which implies that $mis(G_1)=mis(G_2)$. But, if $mis(G'')=\phi(t)$, then by Claim $\ref{cl4.2}$, $R_1\cong C_5$, which implies that it is impossible for $mis(G_1)=mis(G_2)$. This is a contradiction.

If $t$ is odd, then there is no one component with matching number $1$ in $G'-v$, and every non-trivial component in $G'-v$ is isomorphic to $C_5$. Thus,
\[\phi(t)=mis(G')=mis(G'-v)+mis(G'-N_{G'}[v])=5^{\frac{t-1}{2}}+3^{\frac{t-1}{2}}=f(t).\]
Moreover, $G'$ is isomorphic to a graph in $\mathcal{F}_t$.

If $t$ is even, then there is exactly one component with matching number $1$ in $G'-v$, and every component with matching number $2$ in $G'-v$ is isomorphic to $C_5$. Thus,
\[\phi(t)=mis(G')=mis(G'-v)+mis(G'-N_{G'}[v])\leq 2\cdot (5^{\frac{t-2}{2}}+3^{\frac{t-2}{2}})=f(t),\]
which implies that $\phi(t)=f(t)$, and $G'$ is isomorphic to a graph in $\mathcal{F}_t$.  \qed

Now, we have shown that when $t\geq 5$, $\phi(t)=f(t)$, and $G'$ is isomorphic to a graph in $\mathcal{F}_t$. Since $mis(G)=mis(G')=\phi(t)$, we can infer that $G'$, being isomorphic to a graph in $\mathcal{F}_t$, implies that $G$ itself is isomorphic to a graph in $\mathcal{F}_t$.

The proof of Theorem \ref{thm4} is complete. \qed

%%%%%%%%%%%%%%%%%%%%%%%%%%%%%%%%%%%%%%%%%%%%%%%%%%%%%%%%%%%%%%%%%%%%%%

\section*{Use of AI tools declaration}
The authors declare that they have not used Artificial Intelligence (AI) tools in the creation of this article.

\section*{Data Availability}
Data sharing is not applicable to this paper as no datasets were generated or analyzed during
the current study.

\section*{Acknowledgments}
This work was supported by Beijing Natural Science Foundation (No. 1232005). 

\section*{Conflict of interest}
The authors declare that they have no conflicts of interest.

	\bibliographystyle{unsrt}

\begin{thebibliography}{99}
		
        \bibitem{Berge1957}
		C. Berge, Two theorems in graph theory, Proc. Nat. Acad. Sci. U.S.A. 43 (1957), 842--844.
		
		\bibitem{Chang1999}
		G.J. Chang, M.J. Jou, The number of maximal independent sets in connected triangle-free graphs, Discrete Math. 197/198 (1999), 169--178. 
		
		\bibitem{Furedi1987}
		Z. F{\"{u}}redi, The number of maximal independent sets in connected graphs, J. Graph Theory 11 (1987), 463--470.
	
		\bibitem{Griggs1988}
		J.R. Griggs, C.M. Grinstead, D.R. Guichard, The number of maximal independent sets in a connected graph, Discrete Math. 68 (1988), 211--220.
		
    %    \bibitem{Harary1969}
	%	F. Harary, Graph Theory, Addison-Wesley, Reading, MA, 1969.
		
		\bibitem{Hoang2019}		
		D.T. Hoang, T.N. Trung, Coverings, matchings and the number of maximal independent sets of graphs, Australas. J. Combin. 73 (2019), 424--431.
		
		
		\bibitem{Hujtera1993}
		M. Hujtera, Z. Tuza, The number of maximal independent sets in triangle-free graphs, SIAM J. Discrete Math. 6 (1993), 284--288.

		\bibitem{Koh2008}
		K.M. Koh, C.Y. Goh, F.M. Dong, The maximum number of maximal independent sets in unicyclic connected graphs, Discrete Math. 308 (2008), 3761--3769.
		
		\bibitem{Liu1994}
		J. Liu, Constraints on the number of maximal independent sets in graphs, J. Graph Theory 18 (1994), 195--204.
		
		\bibitem{Lovasz1986}
		L. Lov\'{a}sz, M.D. Plummer, Matching Theory, North-Holland, Amsterdam, New York, Oxford, Tokyo, 1986.
		
		\bibitem{Miller1960}
		 R.E. Miller, D.E. Muller, A problem of maximum consistent subsets, IBM Research Report RC-240, 1960, J.T. Watson Research Center, New York, USA.
		
		\bibitem{Moon1965}
		J.W. Moon, L. Moser, On cliques in graphs, Isr. J. Math. 3 (1965), 23--28.
						
		
        \bibitem{Palmer2023}
		C. Palmer, B. Patk\'{o}s, On the number of maximal independent sets: From Moon-Moser to Hujter-Tuza,
		J. Graph Theory 104 (2023), 440--445.

		\bibitem{Sagan1988}
		B.E. Sagan, A note on independent sets in trees, SIAM J. Discrete Math. 1 (1988), 105--108.
		
		\bibitem{Sagan2006}
		B.E. Sagan, V.R. Vatter, Maximal and maximum independent sets in graphs with at most $r$ cycles, J. Graph Theory 53 (2006), 283--314.

		
	   \bibitem{Taleskii2022}
	    D.S. Taletskii, D.S. Malyshev, The number of maximal independent sets in trees with a given number of leaves, Discrete Appl. Math. 314 (2022), 321--330.
		
	
		\bibitem{Wilf1986}
		H.S. Wilf, The number of maximal independent sets in a tree, SIAM J. Alg. Discrete Meth. 7 (1986), 125--130.
		
		\bibitem{Wloch2008}
		I. W{\l}och, Trees with extremal numbers of maximal independent sets including the set of leaves, Discrete Math. 308 (2008), 4768--4772.
		
	
		\bibitem{Ying2006}
		G.C. Ying, K.K. Meng, B.E. Sagan, V.R. Vatter, Maximal independent sets in graphs with at most $r$ cycles, J. Graph Theory 53 (2006), 270--282.
	
	\end{thebibliography}

\end{document}